\documentclass[12pt,a4paper,leqno]{amsart}

\usepackage{anysize}
\usepackage[utf8]{inputenc}
\usepackage{amsmath}
\usepackage{amssymb}
\usepackage{array}

\begin{document}

\title[]{WZ proofs of Ramanujan-type series \\ (via $_2F_1$ evaluations)}
\author{Jesús Guillera}
\address{Department of Mathematics, University of Zaragoza, 50009 Zaragoza, SPAIN}.
\email{}
\date{}

\maketitle

\begin{abstract}
We use Zeilberger's algorithm for proving some identities of Ramanujan-type via $_2F_1$ evaluations.
\end{abstract}

\section{Introduction}
D. Zeilberger wrote the Maple package twoFone, which found automatically many nice formulas, like for example \cite[Theorem 14]{ekhad-zeil}, that with the substitution $n \to -1/4+k$ (valid due to Carlson's theorem) turns into
\[
\sum_{n=0}^{\infty} \frac{\left(\frac14-k\right)_n\left(\frac14-3k \right)_n}{(1+2k)_n(1)_n}(9-4\sqrt{5})^n=\frac{C_1 \cdot 2^{8k}}{5^{2k}(5+2\sqrt{5})^k} \frac{(1)_k\left( \frac12 \right)_k}{\left( \frac{11}{20} \right)_k \left( \frac{19}{20} \right)_k}.
\]
Taking $k=1/4$, we deduce the value of the constant $C_1$:
\[
C_1=\sqrt[4]{\frac{22\sqrt {5} + 50}{125}} \, \frac{\sqrt{\pi}}{ \Gamma\left( \frac{11}{20} \right) \Gamma \left( \frac{19}{20} \right) }.
\]
The main idea of this paper\footnote{For reading this paper is necessary familiarity with Zeilberger's algorithm but only as a black box.} is the following one: Let $A(n,k)$ be the adding term of the above series, and $B(n,k)=A(n,k)(n+bk+c)$, then we get the values of $b$ and $c$ such that the coefficient of $K^2$ is $0$ using the Maple code:
\begin{verbatim}
              coK2:=coeff(Zeilberger(B(n,k),k,n,K)[1],K,2);
              coes:=coeffs(coK2,k);
              solve({coes},{b,c});
\end{verbatim}
For the above example we get $b=(\sqrt{5}-1)/2$ and $c=(5-\sqrt{5})/40$, that leads to the complementary formula
\begin{align}
\sum_{n=0}^{\infty} & \frac{\left(\frac14-k\right)_n\left(\frac14-3k\right)_n}{(1+2k)_n(1)_n}(9-4\sqrt{5})^n \left[ \! \frac{}{} 40n+20(\sqrt{5}-1)k+5-\sqrt{5} \right] \\
& = C_2 \frac{2^{8k}}{5^{2k}(5+2\sqrt{5})^k} \frac{(1)_k\left( \frac12 \right)_k}{\left( \frac{3}{20} \right)_k \left( \frac{7}{20} \right)_k}, 
\end{align}
where we get $C_2$ taking $k=1/4$, and we see that $C_1 C_2=2\sqrt{10+5\sqrt{5}} \, \pi^{-1}$. Substituting $k=0$, and multiplying both series, we obtain
\[
\sum_{n=0}^{\infty} \frac{\left(\frac14\right)_n^2}{(1)_n^2}(9-4\sqrt{5})^n
\sum_{n=0}^{\infty} \frac{\left(\frac14\right)_n^2}{(1)_n^2}(9-4\sqrt{5})^n (40n+5-\sqrt{5})= 
\frac{2\sqrt{10+5\sqrt{5}}}{\pi}.
\]
Finally, using Clausen's identity 
\begin{equation}\label{clausen-id-1} 
\left(\sum_{n=0}^{\infty} \frac{(a)_n(b)_n}{\left(a+b+\frac12\right)_n (1)_n} z^n\right)^2 = 
\sum_{n=0}^{\infty} \frac{(2a)_n(2b)_n (a+b)_n}{\left(a+b+\frac12\right)_n (2a+2b)_n (1)_n} z^n, 
\end{equation}
and its derivative, the product transforms into the Ramanujan-type series
\[
\sum_{n=0}^{\infty} \frac{\left(\frac12\right)_n^3}{(1)_n^3}(9-4\sqrt{5})^n (20n+5-\sqrt{5})=
\frac{2\sqrt{10+5\sqrt{5}}}{\pi}.
\]
We presented the above example in the talk \cite[Slides 29 \& 30]{gui-talk-rutgers} at Rutgers University. In this paper we will use Clausen's formula and also these other known transformations of Clausen type:
\begin{align}
&\left(\sum_{n=0}^{\infty} \frac{(a)_n(b)_n}{\left(\frac{a+b+1}{2}\right)_n (1)_n} z^n\right)^2 = 
\sum_{n=0}^{\infty} \frac{(a)_n(b)_n \left(\frac{a+b}{2} \right)_n}{\left(\frac{a+b+1}{2}\right)_n (a+b)_n (1)_n} (4z(1-z))^n, \label{clausen-id-2} \\
&\left(\sum_{n=0}^{\infty} \frac{(a)_n(b)_n}{(2b)_n (1)_n} z^n\right)^2 = (1-z)^{-a}
\sum_{n=0}^{\infty} \frac{(a)_n(b)_n (2b-a)_n}{\left(b+\frac12\right)_n (2b)_n (1)_n} \left( \frac{z^2}{4(z-1)} \right)^n, \label{clausen-id-3}
\end{align} 
to derive rational identities of the same style. Observe that Clausen's transformations can be proved by Zeilberger's algorithm, thanks to the presence of free parameters. 
\par In \cite{ekhad-zeil} and also in \cite{gessel} there are identities for other values of $z$, and we can use our technique to prove the corresponding Ramanujan \cite{Ramanujan} (or Ramanujan-type) series for $1/\pi$ (see the surveys \cite{BaBeCh} and \cite{Zudilin-rama-wind}, and their references). Related papers with proofs based on Zeilberger's algorithm are \cite{gui-gen-rama, gui-WZ-pairs-rama, gui-proofs-zeil}. Another kind of elementary proofs is shown in \cite{gui-zud-rama-translation}. But, who is the author\footnote{Shalosh Ekhad is the name that Zeilberger gave to his computer.} of \cite{ekhad-zeil}?

\section{Proofs of Ramanujan-type series via $_2F_1$ evaluations}
We will use Zeilberger's algorithm and the Wilf and Zeilberger (WZ) method for proving the formulas of this section.

\subsection{Identities with $s=3$ and $z=1/2$}
\begin{align}
&\sum_{n=0}^{\infty} \frac{\left(\frac12\right)_n\left(\frac13\right)_n\left(\frac23\right)_n}{(1)_n^3} \left(\frac12 \right)^n = \frac{\pi^2 \sqrt[3]{2}(4-2\sqrt 2)} {2 \, \Gamma^2 \!\left(\frac{13}{24}\right)\Gamma^2 \!\left(\frac{19}{24}\right)\Gamma^2 \!\left(\frac23\right)}, \label{Rama-1over2-0} \\
&\sum_{n=0}^{\infty} \frac{\left(\frac12\right)_n\left(\frac13\right)_n\left(\frac23\right)_n}{(1)_n^3} (6n+1) \left(\frac12 \right)^n = \frac{3\sqrt 3}{\pi}. \label{Rama-1over2-1}
\end{align}
\begin{proof}
We can check that Wilf--Zeilberger's algorithm certifies the following identities:
\begin{equation}\label{eq1}
\sum_{n=0}^{\infty} \frac{\left(\frac13+2k\right)_n\left(\frac16+4k\right)_n}{\left(1+3k\right)_n(1)_n } \left( \frac12 \right)^n = C_1\frac{\left(\frac13\right)_k(1)_k}{\left(\frac{13}{24}\right)_k\left(\frac{19}{24}\right)_k} \left(\frac{27}{4}\right)^k, 
\end{equation}
and
\begin{equation}\label{eq2}
\sum_{n=0}^{\infty} \frac{\left(\frac13+2k\right)_n\left(\frac16+4k\right)_n}{\left(1+3k\right)_n(1)_n } \left( \frac12 \right)^n (12n+24k+1)= C_2\frac{\left(\frac13\right)_k(1)_k}{\left(\frac{1}{24}\right)_k\left(\frac{7}{24}\right)_k} \left(\frac{27}{4}\right)^k, 
\end{equation}
where
\[
C_1=\frac{\pi \sqrt[3]{4}\sqrt{4-2\sqrt 2}}{2\Gamma\!\left(\frac{13}{24}\right) \Gamma\!\left(\frac{19}{24}\right) \Gamma \!\left(\frac23\right)}, \qquad C_2=\frac{3\sqrt 3}{C_1 \pi},
\]
where we have determined the constants $C_1$ and $C_2$ by taking $k=-1/6$. Letting $k=0$, squaring (\ref{eq1}) and applying Clausen's identity, we obtain (\ref{Rama-1over2-0}). Letting $k=0$, multiplying (\ref{eq1}) and (\ref{eq2}) and using Clausen's identity (and its derivative), we obtain  (\ref{Rama-1over2-1}).
\end{proof}

\subsection{Identities with $s=4$ and $z=1/9$} 
\begin{align}
&\sum_{n=0}^{\infty} \frac{\left(\frac12\right)_n\left(\frac14\right)_n\left(\frac34\right)_n}{(1)_n^3} \left(\frac19 \right)^n =  (\sqrt{2}+\sqrt{6}+2)^2 \frac{2^\frac{1}{6} \, 3^{\frac54} \, \Gamma^2 \left( \frac78 \right) \Gamma^2 \left( \frac23\right)}{12^2 \, \Gamma^2 \left( \frac58 \right) \Gamma^4 \left( \frac{23}{24} \right)}, \label{Rama-1over9-0} \\
&\sum_{n=0}^{\infty} \frac{\left(\frac12\right)_n\left(\frac14\right)_n\left(\frac34\right)_n}{(1)_n^3} (8n+1) \left(\frac19 \right)^n = \frac{2\sqrt 3}{\pi}. \label{Rama-1over9-1}
\end{align}
\begin{proof}
We can check that Wilf--Zeilberger's algorithm certifies the following identities:
\begin{equation}\label{eq3}
\sum_{n=0}^{\infty} \frac{\left(\frac18-k\right)_n\left(\frac38-k\right)_n}{\left(1+2k\right)_n(1)_n } \left( \frac19 \right)^n = C_1\frac{\left(\frac12\right)_k(1)_k}{\left(\frac{13}{24}\right)_k\left(\frac{23}{24}\right)_k} \left(\frac{2^8}{3^5}\right)^k, 
\end{equation}
and
\begin{equation}\label{eq4}
\sum_{n=0}^{\infty} \frac{\left(\frac18-k\right)_n\left(\frac38-k\right)_n}{\left(1+2k\right)_n(1)_n } \left( \frac19 \right)^n (16n+8k+1)= C_2\frac{\left(\frac12\right)_k(1)_k}{\left(\frac{5}{24}\right)_k\left(\frac{7}{24}\right)_k} \left(\frac{2^8}{3^5}\right)^k, 
\end{equation}
where 
\[
C_1=  (\sqrt{2}+\sqrt{6}+2) \frac{2^\frac{1}{12} \, 3^{\frac58} \, \Gamma \left( \frac78 \right) \Gamma\left( \frac23\right)}{12 \, \Gamma \left( \frac58 \right) \Gamma^2 \left( \frac{23}{24} \right)}, \qquad
C_2=\frac{2\sqrt 3}{C_1 \pi},
\]
have been determined by taking $k=1/8$. Letting $k=0$, squaring (\ref{eq3}) and applying Clausen's identity, we obtain (\ref{Rama-1over9-0}). Letting $k=0$, multiplying (\ref{eq3}) and (\ref{eq4}) and using Clausen's identity (and its derivative), we obtain  (\ref{Rama-1over9-1}).
\end{proof}

\subsection{Identities with $s=4$ and $z=32/81$} 
\begin{align}
&\sum_{n=0}^{\infty} \frac{\left(\frac12\right)_n\left(\frac14\right)_n\left(\frac34\right)_n}{(1)_n^3} \left(\frac{32}{81} \right)^n =  \frac{\sqrt{3} \, \pi}{2 \, \Gamma^2\left( \frac{7}{12}\right) \Gamma^2\left( \frac{11}{12} \right)}, \label{Rama-32over81-0} \\
&\sum_{n=0}^{\infty} \frac{\left(\frac12\right)_n\left(\frac14\right)_n\left(\frac34\right)_n}{(1)_n^3}  (7n+1) \left(\frac{32}{81} \right)^n= \frac{9}{2\pi}. \label{Rama-32over81-1}
\end{align}

\begin{proof}
We can check that Wilf--Zeilberger's algorithm certifies the following identities:
\begin{equation}\label{eq5}
\sum_{n=0}^{\infty} \frac{\left(\frac14-k\right)_n\left(\frac34-k\right)_n}{\left(1+2k\right)_n(1)_n } \left( \frac19 \right)^n = C_1\frac{\left(\frac12\right)_k(1)_k}{\left(\frac{7}{12}\right)_k\left(\frac{11}{12}\right)_k} \left(\frac{2^8}{3^5}\right)^k, \quad C_1= \frac{\sqrt[4]{3} \, \sqrt{2\pi}}{2 \, \Gamma\left( \frac{7}{12}\right) \Gamma\left( \frac{11}{12} \right)}, 
\end{equation}
and
\begin{equation}\label{eq6}
\sum_{n=0}^{\infty} \frac{\left(\frac14-k\right)_n\left(\frac34-k\right)_n}{\left(1+2k\right)_n(1)_n } \left( \frac19 \right)^n (16n+20k+1)= C_2\frac{\left(\frac12\right)_k(1)_k}{\left(\frac{1}{12}\right)_k\left(\frac{5}{12}\right)_k} \left(\frac{2^8}{3^5}\right)^k, \quad C_2=\frac{9}{2 C_1 \pi},
\end{equation}
where we have determined the constants $C_1$ and $C_2$ by taking $k=1/4$. Letting $k=0$, squaring (\ref{eq5}) and applying Clausen's (\ref{clausen-id-2}) identity, we obtain (\ref{Rama-32over81-0}). Letting $k=0$, multiplying (\ref{eq5}) and (\ref{eq6}), and finally applying Clausen's (\ref{clausen-id-2}) identity (and its derivative), we obtain  (\ref{Rama-32over81-1}).
\end{proof}

\subsection{Identities with $s=3$ and $z=-9/16$} 
\begin{align}
&\sum_{n=0}^{\infty} \frac{\left(\frac12\right)_n\left(\frac13\right)_n\left(\frac23\right)_n}{(1)_n^3} \left(-\frac{9}{16} \right)^n = \frac{16 \pi^2}{27 \, \Gamma^6\left( \frac23 \right)}, \label{Rama--1over8-0} \\
&\sum_{n=0}^{\infty} \frac{\left(\frac12\right)_n\left(\frac13\right)_n\left(\frac23\right)_n}{(1)_n^3}  (5n+1) \left(-\frac{9}{16} \right)^n= \frac{4\sqrt{3}}{3\pi}. \label{Rama--1over8-1}
\end{align}

\begin{proof}
We can check that Wilf--Zeilberger's algorithm certifies the following identities:
\begin{equation}\label{eq7}
\sum_{n=0}^{\infty} \frac{\left(\frac13-k\right)_n\left(\frac23-k\right)_n}{\left(1+k\right)_n(1)_n } \left( -\frac18 \right)^n = C_1\frac{(1)_k}{\left(\frac{5}{6}\right)_k} \left(\frac{3^3}{2^5} \right)^k, \quad C_1= \frac{4 \pi \sqrt{3}}{9 \, \Gamma^3 \left( \frac23 \right)}, 
\end{equation}
and
\begin{equation}\label{eq8}
\sum_{n=0}^{\infty} \frac{\left(\frac13-k\right)_n\left(\frac23-k\right)_n}{\left(1+k\right)_n(1)_n} (9n+3k+1) \left( -\frac18 \right)^n = C_2\frac{(1)_k}{\left(\frac{1}{6}\right)_k} \left(\frac{3^3}{2^5}\right)^k, \quad C_2=\frac{3 \, \Gamma^3 \left( \frac23\right)}{\pi^2},
\end{equation}
where we have determined the constants $C_1$ and $C_2$ by taking $k=1/3$. Letting $k=0$, squaring (\ref{eq7}) and applying Clausen's (\ref{clausen-id-2}) identity, we obtain (\ref{Rama--1over8-0}). Letting $k=0$, multiplying (\ref{eq7}) and (\ref{eq8}), and finally applying Clausen's (\ref{clausen-id-2}) identity (and its derivative), we obtain  (\ref{Rama--1over8-1}).
\end{proof}

\subsection{Identities with $s=2$ and $z=-1$}
\begin{equation}\label{Rama-n1}
\sum_{n=0}^{\infty} \frac{\left(\frac12\right)_n^3}{(1)_n^3} (-1)^n = \frac{2+\sqrt{2}}{4} \, \frac{\Gamma^2 \left(\frac34 \right)}{\Gamma^4 \left(\frac78 \right)}, \qquad
\sum_{n=0}^{\infty} \frac{\left(\frac12\right)_n^3}{(1)_n^3}  (4n+1) (-1)^n = \frac{2}{\pi}. 
\end{equation}
\begin{proof}
We can check that Wilf--Zeilberger's algorithm certifies the following identities:
\begin{equation}\label{eq-n1-1}
\sum_{n=0}^{\infty} \frac{\left(\frac12\right)_n\left(\frac12+2k\right)_n\left(\frac12+4k\right)_n}{(1)_n (1+2k)_n (1+4k)_n } (-1)^n = C \, \frac{(1)_k^2(\frac12)_k^2}{\left(\frac58\right)_k^2 \left( \frac78\right)_k^2}, \quad C=\frac{2+\sqrt{2}}{4} \, \frac{\Gamma^2 \left(\frac34 \right)}{\Gamma^4 \left(\frac78 \right)},
\end{equation}
and
\begin{equation}\label{eq-n1-2}
\sum_{n=0}^{\infty} \frac{\left(\frac12\right)_n\left(\frac12+k\right)_n\left(\frac12+2k\right)_n}{\left(1\right)_n (1+k)_n (1+2k)_n } (-1)^n (4n+4k+1) =
\frac{2}{\pi} \, \frac{(1)_k^2}{\left(\frac14\right)_k\left(\frac34\right)_k}.
\end{equation}
Taking $k=0$ in (\ref{eq-n1-1} and (\ref{eq-n1-2}), we obtain the identities in (\ref{Rama-n1}).
\end{proof}

\subsection{Identities with $s=2$ and $z=1/4$}
\begin{equation}\label{Rama-1over4}
\sum_{n=0}^{\infty} \frac{\left(\frac12\right)_n^3}{(1)_n^3} \left(\frac14 \right)^n = \frac{8 \, \pi^2 \sqrt[3]{2} \, \sqrt{3}}{27 \, \Gamma^6 \left( \frac23\right)}, \qquad
\sum_{n=0}^{\infty} \frac{\left(\frac12\right)_n^3}{(1)_n^3}  (6n+1) \left(\frac14 \right)^n = \frac{4}{\pi}. 
\end{equation}
\begin{proof}
We can check that Wilf--Zeilberger's algorithm certifies the following identities:
\begin{equation}\label{eq9}
\sum_{n=0}^{\infty} \frac{\left(\frac12-k\right)_n\left(\frac12+3k\right)_n\left(\frac12+k\right)_n}{\left(1+k\right)_n (1+2k)_n (1)_n } \left( \frac14 \right)^n = C \, \frac{(1)_k^2}{(\frac56)_k^2} \left( \frac{16}{27}\right)^k, \quad C=\frac{8 \, \pi^2 \sqrt[3]{2} \, \sqrt{3}}{27 \, \Gamma^6 \left( \frac23\right)},
\end{equation}
and
\begin{equation}\label{eq10}
\sum_{n=0}^{\infty} \frac{\left(\frac12-k\right)_n\left(\frac12+3k\right)_n\left(\frac12+k\right)_n}{\left(1+k\right)_n (1+2k)_n (1)_n } \left( \frac14 \right)^n (6n+6k+1) =
\frac{4}{\pi} \, \left(\frac{16}{27}\right)^k \frac{(1)_k^2}{\left(\frac16\right)_k\left(\frac56\right)_k}.
\end{equation}
Taking $k=0$ in (\ref{eq9}) and (\ref{eq10}), we obtain the identities (\ref{Rama-1over4}).
\end{proof}

\subsection{Identities with $s=2$ and $z=-1/8$}
\begin{equation}\label{Rama-n1over8}
\sum_{n=0}^{\infty} \frac{\left(\frac12\right)_n^3}{(1)_n^3} \left(-\frac18 \right)^n = \frac{\pi \, \sqrt{6}}{3 \, \Gamma^2 \left( \frac{7}{12} \right)\Gamma^2 \left( \frac{11}{12} \right)}, \qquad
\sum_{n=0}^{\infty} \frac{\left(\frac12\right)_n^3}{(1)_n^3}  (6n+1) \left(-\frac18 \right)^n= \frac{2\sqrt{2}}{\pi}. 
\end{equation}
\begin{proof}
We can check that Wilf--Zeilberger's algorithm certifies the following identities:
\begin{equation}\label{eq11}
\sum_{n=0}^{\infty} \frac{\left(\frac12-2k\right)_n\left(\frac12+6k\right)_n\left(\frac12+2k\right)_n}{\left(1+2k\right)_n (1+4k)_n (1)_n } \left( -\frac18 \right)^n = C_1 \, \frac{(\frac12)_k^2(1)_k^2}{(\frac{7}{12})_k^2(\frac{11}{12})_k^2} \left( \frac{32}{27}\right)^{2k},  
\end{equation}
where
\[
C_1=\frac{\pi \, \sqrt{6}}{3 \, \Gamma^2 \left( \frac{7}{12} \right)\Gamma^2 \left( \frac{11}{12} \right)},
\]
and
\begin{equation}\label{eq12}
\sum_{n=0}^{\infty} \frac{\left(\frac12-k\right)_n\left(\frac12+3k\right)_n\left(\frac12+k\right)_n}{\left(1+k\right)_n (1+2k)_n (1)_n } \left( -\frac18 \right)^n (6n+6k+1) =
\frac{2\sqrt{2}}{\pi} \, \left(\frac{32}{27}\right)^k \frac{(1)_k^2}{\left(\frac16\right)_k\left(\frac56\right)_k}.
\end{equation}
Taking $k=0$ in (\ref{eq11}) and (\ref{eq11}), we obtain the identities (\ref{Rama-n1over8}).
\end{proof}

\subsection{Identities with $s=4$ and $z=-1/4$}
\begin{align}\label{Rama-n1over4}
&\sum_{n=0}^{\infty} \frac{\left(\frac12\right)_n\left(\frac14\right)_n\left(\frac34\right)_n}{(1)_n^3} \left(-\frac14 \right)^n =\frac{2}{25} \frac{\sqrt[4]{5} \, (5+\sqrt{5}) \, \pi}{\Gamma^2\left( \frac{11}{20}\right)\Gamma^2\left( \frac{19}{20}\right)}, \\
&\sum_{n=0}^{\infty} \frac{\left(\frac12\right)_n\left(\frac14\right)_n\left(\frac34\right)_n}{(1)_n^3}  (20n+3) \left(-\frac14 \right)^n = \frac{8}{\pi}. 
\end{align}
\begin{proof}
We can check that Wilf--Zeilberger's algorithm certifies the following identities:
\begin{equation}\label{eq13}
\sum_{n=0}^{\infty} \frac{\left(\frac12+2k\right)_n\left(\frac14-k\right)_n\left(\frac34+5k\right)_n}{\left(1+2k\right)_n (1+4k)_n (1)_n } \left( -\frac14 \right)^n = C \, \frac{(\frac12)_k^2(1)_k^2}{(\frac{11}{20})_k^2(\frac{19}{20})_k^2} \left( \frac{2^{12}}{5^5}\right)^k, 
\end{equation}
where
\[
C=\frac{2}{25} \frac{\sqrt[4]{5} \, (5+\sqrt{5}) \, \pi}{\Gamma^2\left( \frac{11}{20}\right)\Gamma^2\left( \frac{19}{20}\right)},
\]
and
\begin{multline}\label{eq14}
\sum_{n=0}^{\infty} \frac{\left(\frac12+2k\right)_n\left(\frac14-k\right)_n\left(\frac34+5k\right)_n}{\left(1+2k\right)_n (1+4k)_n (1)_n } \left( -\frac14 \right)^n (20n+20k+3) \\ =
\frac{8}{\pi} \, \left(\frac{2^{12}}{5^5}\right)^k \frac{(1)_k^2\left(\frac12 \right)_k^2}{\left(\frac{3}{20}\right)_k\left(\frac{7}{20}\right)_k\left(\frac{11}{20}\right)_k\left(\frac{19}{20}\right)_k}.
\end{multline}
Taking $k=0$ in (\ref{eq13}) and (\ref{eq14}), we obtain the identities (\ref{Rama-n1over4}).
\end{proof}

\subsection{Identities with $s=3$ and $z=-4$}
\begin{align}\label{Rama-n4}
& \sum_{n=0}^{\infty} \frac{\left(\frac12\right)_n\left(\frac13\right)_n\left(\frac23\right)_n}{(1)_n^3} (-4)^n  \, \,``\!=\!" \, \, \frac{\sqrt[3]{5} \, \sqrt[5]{81}}{5} \, \left( \frac{\Gamma \left(\frac23 \right)\Gamma \left(\frac35 \right)\Gamma \left(\frac45 \right)}{\Gamma \left(\frac{8}{15} \right) \Gamma \left(\frac{11}{15} \right) \Gamma \left(\frac{13}{15} \right) \Gamma \left(\frac{14}{15} \right)} \right)^2, \\
& \sum_{n=0}^{\infty} \frac{\left(\frac12\right)_n \left(\frac13\right)_n\left(\frac23\right)_n}{(1)_n^3}  (4n+1)  (-4)^n  \, \, ``\!=\!" \, \, \frac{3 \sqrt{3}}{\pi}. \nonumber
\end{align}
Here $``\!=\!"$ means the analytic continuation of the corresponding function of $z$ at $z=-4$.
\begin{proof}
We can check that Wilf--Zeilberger's algorithm certifies the following identities:
\begin{equation}\label{eq-n4-1}
\sum_{n=0}^{\infty} \frac{\left(\frac12+3k\right)_n\left(\frac13+k\right)_n\left(\frac23+5k\right)_n}{(1+3k)_n (1+6k)_n (1)_n } (-4)^n \, \,``\!=\!" \, \, C \, \frac{(1)_k^2(\frac23)_k^2}{\left(\frac{11}{15}\right)_k^2 \left( \frac{14}{15}\right)_k^2} \left( \frac{3^6}{5^5} \right)^k,
\end{equation}
where
\[
C=\frac{\sqrt[3]{5} \, \sqrt[5]{81}}{5} \, \left( \frac{\Gamma \left(\frac23 \right)\Gamma \left(\frac35 \right)\Gamma \left(\frac45 \right)}{\Gamma \left(\frac{8}{15} \right) \Gamma \left(\frac{11}{15} \right) \Gamma \left(\frac{13}{15} \right) \Gamma \left(\frac{14}{15} \right)} \right)^2,
\]
is determined by taking $k=-2/15$, and
\begin{multline}\label{eq-n4-2}
\sum_{n=0}^{\infty} \frac{\left(\frac12+3k\right)_n\left(\frac13+k\right)_n\left(\frac23+5k\right)_n}{\left(1+3k\right)_n (1+6k)_n (1)_n } (-4)^n (15n+30k+4) \\ \, \,``\!=\!" \, \,
\frac{3 \sqrt{3}}{\pi} \, \frac{(1)_k^2 \left( \frac23 \right)_k^2}{\left(\frac{2}{15} \right)_k \left(\frac{8}{15} \right)_k \left(\frac{11}{15} \right)_k \left(\frac{14}{15} \right)_k} \left( \frac{3^6}{5^5} \right)^k.
\end{multline}
Taking $k=0$ in (\ref{eq-n4-1}) and (\ref{eq-n4-2}), we obtain the identities in (\ref{Rama-n4}).
\end{proof}

\subsection{Identities with $s=3$ and $z=2/27$}
\begin{align}
&\sum_{n=0}^{\infty} \frac{\left(\frac12\right)_n\left(\frac13\right)_n\left(\frac23\right)_n}{(1)_n^3} \left(\frac{2}{27}\right)^n =\frac{3\pi}{4 \, \Gamma^2 \left(\frac23\right)\Gamma^2 \left( \frac56\right)},  \label{Rama-2-27} \\
&\sum_{n=0}^{\infty} \frac{\left(\frac12\right)_n\left(\frac13\right)_n\left(\frac23\right)_n}{(1)_n^3} (15n+2) \left(\frac{2}{27}\right)^n = \frac{27}{4\pi}, \nonumber
\end{align}

The proof is inspired by the formulas (30.4a) and (30.6a) of \cite{gessel} with $z=2/27$, by making the substitutions $n \to -1/6+k$ and $n \to -5/6+k$ respectively. 
\begin{proof}[Proof of (\ref{Rama-2-27})]
We can check directly that Wilf--Zeilberger's algorithm certifies the following identities:
\begin{equation}\label{eq15}
\sum_{n=0}^{\infty} \frac{\left(\frac23\right)_n\left(\frac56-k\right)_n}{\left(1+k\right)_n(1)_n } \frac{\left(\frac13+2k\right)_n}{\left(\frac13\right)_n} \frac{15n+3k+5}{3n+1}\left( \frac{2}{27} \right)^n = \frac{\left(1\right)_k}{\left(\frac23\right)_k} \, \frac92 \, \frac{\sqrt{\pi}}{\Gamma \left(\frac23 \right)\Gamma \left(\frac56 \right)},
\end{equation}
and
\begin{equation}\label{eq16}
\sum_{n=0}^{\infty} \frac{\left(\frac13\right)_n\left(\frac16-k\right)_n}{\left(1+k\right)_n(1)_n } \frac{\left(\frac23+2k\right)_n}{\left(\frac23\right)_n} (15n+3k+1) \left( \frac{2}{27} \right)^n = \frac{\left(1\right)_k}{\left(\frac13\right)_k} \, \frac{9\sqrt{3}}{4 \pi} \, \frac{\Gamma \left(\frac23 \right)\Gamma \left(\frac56 \right)}{\sqrt{\pi}}.
\end{equation}
Taking $k=0$ in (\ref{eq15}), using Euler's transformation:
\begin{equation}\label{euler-id} 
\sum_{n=0}^{\infty} \frac{(a)_n(b)_n}{(c)_n(1)_n} z^n = (1-z)^{-a} \sum_{n=0}^{\infty} \frac{(a)_n(c-b)_n}{(c)_n(1)_n} \left( \frac{z}{z-1}\right)^n, 
\end{equation}
and finally applying Clausen's identity (\ref{clausen-id-1}), we obtain the first formula in (\ref{Rama-2-27}). Taking $k=0$, multiplying (\ref{eq15}) and (\ref{eq16}), using Euler's transformation and applying Clausen's formula (\ref{clausen-id-1}) and its derivative, we obtain the second formula in (\ref{Rama-2-27}).
\end{proof}

\subsection{Identities with $s=6$ and $z=-9/40^3$}
\begin{align}
&\sum_{n=0}^{\infty} \frac{\left(\frac12\right)_n\left(\frac16\right)_n\left(\frac56\right)_n}{(1)_n^3} \left(-\frac{9}{40^3} \right)^n = \frac{8 \, \sqrt[6]{2} \, \sqrt{5} \, \pi}{27 \, \Gamma^2\left(\frac23\right) \Gamma^2\left(\frac56\right)}, \label{Rama-9-64000} \\
&\sum_{n=0}^{\infty} \frac{\left(\frac12\right)_n\left(\frac16\right)_n\left(\frac56\right)_n}{(1)_n^3} (4554n+279) \left(-\frac{9}{40^3}\right)^n = \frac{160\sqrt{30}}{\pi}, \nonumber
\end{align}
\begin{proof}
We can check directly that Wilf--Zeilberger's algorithm certifies the following identities:
\begin{equation}\label{eq17}
\sum_{n=0}^{\infty} \frac{\left(\frac12-k \right)_n\left(\frac56+k\right)_n}{\left(1+\frac{k}{2} \right)_n(1)_n } \frac{\left(\frac12+k\right)_n}{\left(\frac12+\frac{k}{2} \right)_n} \frac{10n+6k+5}{10n+5k+5}\left( \frac{3}{128} \right)^n = C_1 \, \frac{(1)_k}{\left( \frac56 \right)_k} \left(\frac89\right)^k,
\end{equation}
and
\begin{multline}\label{eq18}
\sum_{n=0}^{\infty} \frac{\left(\frac12-k\right)_n\left(\frac16+k\right)_n}{\left(1+\frac{k}{2}\right)_n(1)_n } \frac{\left(\frac12+k\right)_n}{\left(\frac12+\frac{k}{2}\right)_n} \left( \frac{3}{128} \right)^n \\ \times \frac{100n^2+56n+3+k(80n+12k+20)}{2n+k+1}   =  C_2 \, \frac{(1)_k}{\left( \frac16 \right)_k} \left(\frac89\right)^k.
\end{multline}
where
\[
C_1=\frac{16}{45} \, \frac{\sqrt{6\pi}}{\Gamma\left(\frac23\right)\Gamma\left(\frac56\right)}, \qquad 
C_2=\frac{24\sqrt{2} \, \Gamma\left(\frac23\right)\Gamma\left(\frac56\right)}{3 \pi \sqrt{\pi}}.
\]
Taking $k=0$ in (\ref{eq17}), using Euler's transformation (\ref{euler-id}), and finally applying Clausen's identity (\ref{clausen-id-3}), we obtain the first formula in (\ref{Rama-9-64000}). Taking $k=0$, multiplying (\ref{eq17}) and (\ref{eq18}), using Euler's transformation and applying Clausen's formula (\ref{clausen-id-3}) and its derivative, we obtain the second formula in (\ref{Rama-9-64000}).
\end{proof}

\subsection{Identities with $s=4$ and $z=-1/48$}
\begin{align}
&\sum_{n=0}^{\infty} \frac{\left(\frac12\right)_n\left(\frac14\right)_n\left(\frac34\right)_n}{(1)_n^3} \left(-\frac{1}{48} \right)^n = \frac{2 \, \sqrt{2} \, \pi}{3\,\sqrt[4]{3} \, \Gamma^4\left( \frac34\right)}, \label{Rama-1-48} \\
&\sum_{n=0}^{\infty} \frac{\left(\frac12\right)_n\left(\frac14\right)_n\left(\frac34\right)_n}{(1)_n^3} (28n+3) \left(-\frac{1}{48}\right)^n = \frac{16\sqrt{3}}{3\pi}, \nonumber
\end{align}
\begin{proof}
We can check directly that Wilf--Zeilberger's algorithm certifies the following identities:
\begin{equation}\label{eq19}
\sum_{n=0}^{\infty} \frac{\left(\frac12+2k\right)_n\left(\frac14-k\right)_n}{\left(1+2k\right)_n(1)_n} \left(-\frac13 \right)^n = C_1 \frac{\left(\frac12\right)_k(1)_k}{\left(\frac{7}{12}\right)_k\left(\frac{11}{12}\right)_k} \left(\frac{4}{3}\right)^k, 
\end{equation}
and
\begin{equation}\label{eq20}
\sum_{n=0}^{\infty} \frac{\left(\frac12+2k \right)_n\left(\frac14-k\right)_n}{\left(1+2k\right)_n(1)_n} \left( -\frac13 \right)^n (8n+4k+1)= C_2\frac{\left(\frac12\right)_k(1)_k}{\left(\frac{1}{4}\right)_k^2} \left(\frac{4}{3}\right)^k, 
\end{equation}
where
\[
C_1= \frac{\sqrt{6\pi}}{3\, \Gamma^2\left( \frac34\right)}, \qquad C_2=\frac{2 \, \sqrt[4]{3} \, \Gamma^2\left( \frac34\right)}{\pi \sqrt{\pi}}.
\]
Taking $k=0$ in (\ref{eq19}), using Euler's transformation (\ref{euler-id}), and finally applying Clausen's identity (\ref{clausen-id-3}), we obtain the first formula in (\ref{Rama-1-48}). Taking $k=0$, multiplying (\ref{eq19}) and (\ref{eq20}), using Euler's transformation and applying Clausen's formula (\ref{clausen-id-3}) and its derivative, we obtain the second formula in (\ref{Rama-1-48}).
\end{proof}

\section*{APPENDIX}

\subsection*{A rapid formula for the lemniscate constant $\omega$}

\begin{equation}
\frac{1} {24\sqrt6 \sqrt[8]{125}}\sum_{n=0}^{\infty} \frac{(\frac18)_n(\frac38)_n}{(1)_n^2} (41+1288n) \left(\frac{-1}{5\cdot 72^2}\right)^n = \frac{1}{\omega}.
\end{equation}

\subsection*{Proof of the formula}

Using Clausen's identity, we get
\begin{equation}\label{id1}
\sum_{n=0}^{\infty} \frac{(\frac12)_n(\frac14)_n(\frac34)_n}{(1)_n^3} z^n=
\left(\sum_{n=0}^{\infty} \frac{(\frac18)_n(\frac38)_n}{(1)_n^2} z^n \right)^2.
\end{equation}
Applying the operator $z \, d/dz$, we obtain
\begin{equation}\label{id2}
\sum_{n=0}^{\infty} \frac{(\frac12)_n(\frac14)_n(\frac34)_n}{(1)_n^3} n z^n=
2 \sum_{n=0}^{\infty} \frac{(\frac18)_n(\frac38)_n}{(1)_n^2} n z^n   \sum_{n=0}^{\infty} \frac{(\frac18)_n(\frac38)_n}{(1)_n^2} z^n.
\end{equation}
From (\ref{id1}) and (\ref{id2}), we see that
\begin{multline}
a \sum_{n=0}^{\infty} (\frac{(\frac12)_n(\frac14)_n(\frac34)_n}{(1)_n^3} z^n + b \sum_{n=0}^{\infty} \frac{(\frac12)_n(\frac14)_n(\frac34)_n}{(1)_n^3} n z^n \\ = \sum_{n=0}^{\infty} \frac{(\frac18)_n(\frac38)_n}{(1)_n^2} z^n \left( a \sum_{n=0}^{\infty} \frac{(\frac18)_n(\frac38)_n}{(1)_n^2} z^n + b \sum_{n=0}^{\infty} \frac{(\frac18)_n(\frac38)_n}{(1)_n^2} 2 n z^n \right).
\end{multline}
We obtain,
\begin{equation}
\sum_{n=0}^{\infty} \frac{(\frac12)_n(\frac14)_n(\frac34)_n}{(1)_n^3} (a+bn) z^n \\ = \sum_{n=0}^{\infty} \frac{(\frac18)_n(\frac38)_n}{(1)_n^2} z^n \left( \sum_{n=0}^{\infty} \frac{(\frac18)_n(\frac38)_n}{(1)_n^2} (a+2bn) z^n \right).
\end{equation}
The formula
\begin{equation}
\sum_{n=0}^{\infty} \frac{(\frac12)_n(\frac14)_n(\frac34)_n}{(1)_n^3} (41+644n) \left(\frac{-1}{5\cdot 72^2}\right)^n = \frac{288 \sqrt{5}}{5\pi}.
\end{equation}
is a well known Ramanujan-type series for $1/\pi$ \cite[see tables]{Gui-meth-rama}. Hence
\begin{equation}
\sum_{n=0}^{\infty} \frac{(\frac18)_n(\frac38)_n}{(1)_n^2} \left(\frac{-1}{5\cdot 72^2}\right)^n \sum_{n=0}^{\infty} \frac{(\frac18)_n(\frac38)_n}{(1)_n^2} (41+1288n) \left(\frac{-1}{5\cdot 72^2}\right)^n=\frac{288 \sqrt{5}}{5\pi}.
\end{equation}
Applying Euler transformation 
\[
{}_2F_1\biggl(\begin{matrix}
a, & b \\
& c \end{matrix} \biggm| z \biggr)=
(1-z)^{-a} 
{}_2F_1\biggl(\begin{matrix}
a, & c-b \\
& c \end{matrix} \biggm| \frac{z}{z-1} \biggr)
\]
to the formula \cite[($A'''$ .1) of Table 4]{Ebisu-hyper-2F1}, we have
\begin{equation}
\sum_{n=0}^{\infty} \frac{(\frac18)_n(\frac38)_n}{(1)_n^2} \left(\frac{-1}{5\cdot 72^2}\right)^n = \frac{2\sqrt3 \sqrt[8]{5} \sqrt{\pi}}{5\Gamma(\frac34)^2},
\end{equation}
and finally, we obtain
\begin{equation}
\sum_{n=0}^{\infty} \frac{(\frac18)_n(\frac38)_n}{(1)_n^2} (41+1288n) \left(\frac{-1}{5\cdot 72^2}\right)^n = \frac{48\sqrt3 \sqrt[8]{5^3} \Gamma(\frac34)^2}{\pi^\frac32} = \frac{24\sqrt6 \sqrt[8]{125}}{\omega},
\end{equation}
where $\omega$ is the lemniscate constant. The above formula is suitable to use with y-cruncher \cite{y-cruncher}.

\end{document}